# Mathématiques occidentales en Chine du XVI[e] au XX[e] siècle[1]

Rémi Anicotte[2]

**Résumé :** Du XVI[e] au XVIII[e] siècle, les jésuites sous patronage portugais puis français ont assuré leur influence en Chine en s'appuyant sur les sciences européennes et en participant à une réforme du calendrier chinois. Leurs traductions d'ouvrages européens incitèrent à développer du lexique chinois de la géométrie et provoquèrent un regain d'intérêt des Chinois pour certains aspects de leur tradition mathématique. Au XIX[e] siècle, les calculs algébrique et différentiel furent introduits par des missionnaires protestants anglais et Li Shanlan inventa sa propre symbolique pour les présenter en chinois. Finalement les notations mathématiques occidentales furent adoptées au XX[e] siècle.

**Mots clés :** Notations mathématiques en Chine ; Numération en Chine ; Lexique chinois des mathématiques ; Lexique japonais des mathématiques ; Calendrier chinois ; Traduction de textes mathématiques; Jésuites en Chine ; Missionnaires protestants en Chine ; Matteo Ricci ; Clavius ; Li Shanlan.

---

Des archéologues ont retrouvé en Chine des manuscrits de mathématiques datant du II[e] siècle avant notre ère[3]. Ils contiennent des procédures de calculs avec des entiers et des fractions, des calculs de longueurs, d'aires et de volumes, des résolutions de problèmes concernant des productions agricoles et artisanales et leurs taxations. Les calculs étaient effectués à l'aide de bâtonnets de calculs qui représentaient les chiffres sur une surface plane, il s'agissait d'un système décimal et positionnel.

Chiffres avec des bâtonnets de calcul

|  | 1 | 2 | 3 | 4 | 5 | 6 | 7 | 8 | 9 |
|---|---|---|---|---|---|---|---|---|---|
| Chiffres utilisés pour les unités, les centaines, les dizaines de milliers, etc. | \| | \|\| | \|\|\| | \|\|\|\| | \|\|\|\|\| | ⊤ | ⊤\| | ⊤\|\| | ⊤\|\|\| |
| Chiffres utilisés pour les dizaines, les milliers, les centaines de milliers, etc. | — | = | ≡ | ≣ | ≣ | ⊥ | ⊥ | ⊥ | ⊥ |

Cette culture mathématique antique fut transmise aux générations suivantes par les *Neuf Chapitres*[4], un ouvrage compilé sous le règne de la dynastie des Han orientaux (23-220). Les mathématiques en Chine continuèrent ensuite à évoluer avec notamment l'apparition de nouvelles méthodes de résolution d'équations au XIII[e] siècle et le développement du calcul sur boulier au XIV[e] siècle.

Chiffres sur un boulier

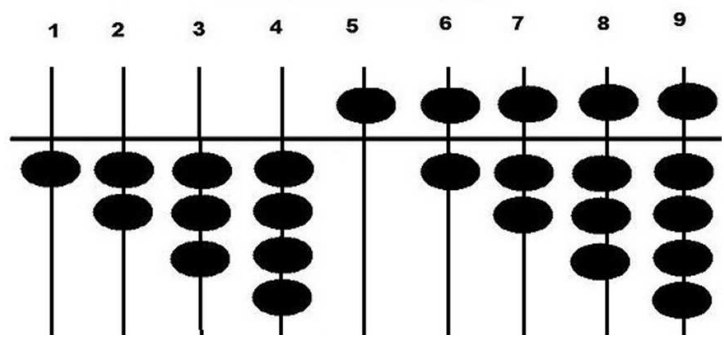

Puis le XVI[e] siècle a vu l'intensification des échanges mondiaux grâce aux navigateurs portugais, en même temps que l'accélération du développement scientifique en Europe. Les jésuites, membres de l'ordre religieux fondé par Ignace de Loyola (1491-1556), arrivèrent en Chine dans ce contexte, d'abord sous patronage portugais, puis à partir de la fin du XVII[e] siècle aussi sous patronage français. Dans les pays catholiques, les jésuites dirigeaient des établissements scolaires qui accordaient une place importante aux matières scientifiques. Signalons, par exemple, Christophorus Clavius (1538-1612), un père jésuite qui participa à la conception du calendrier grégorien que nous utilisons toujours aujourd'hui. Il fut l'auteur, entre autres, d'une arithmétique pratique (*Epitome arithmeticae practicae*) et d'une édition latine des *Éléments* d'Euclide ; ces deux ouvrages ne tardèrent pas à être traduits en chinois.

---

[3] Notamment le *Suan shu shu* découvert dans le Hubei pendant l'hiver 1983-84.
[4] Ce classique des mathématiques en Chine est accessible dans Chemla & Guo (2004).

En Chine justement, la première mission jésuite était établie dans la colonie portugaise de Macao. Le père Matteo Ricci[5] (1522-1610) y arriva en 1582 ; au Collège romain il avait été l'élève de Clavius. Ricci se fit connaître du public chinois grâce à une traduction de la carte du monde connu des européens et à la publication d'un opuscule sur la méthode de mémorisation du *palais de mémoire*[6]. Il bâtit ensuite le projet de gagner l'adhésion des lettrés chinois au christianisme en s'appuyant sur l'astronomie européenne. Il s'agissait de prévoir l'apparition des éclipses et d'améliorer la façon de compenser le décalage entre l'année calendaire chinoise et l'année solaire en optimisant la durée et les dates d'insertion des mois intercalaires. Et c'est en qualité de conseiller sur ces questions que Matteo Ricci fut appelé à Pékin auprès de l'administration impériale en 1601. Après maintes péripéties, c'est le père Johann Adam Schall von Bell (1591-1666) qui mena à terme la conception du calendrier Chongzhen (du nom du dernier empereur Ming) et à la publication du *Chongzhen lishu*, une collection d'écrits d'astronomie calendaire, finalement en 1644, Shunzhi (le premier empereur de la dynastie Qing) promulgua le nouveau calendrier et nomma Schall von Bell à la tête du Bureau impérial d'astronomie.

Hors du champ des mathématiques calendaires, Matteo Ricci, Xu Guangqi (1562-1633) et Li Zhizao (1565-1630) traduisirent en chinois des ouvrages de Clavius. En 1607 sortirent les six premiers livres des *Éléments* d'Euclide, ceux consacrés à la géométrie et aux rapports. Mais cette première traduction n'explicitait pas ce que sont les définitions, axiomes, propositions et preuves. Elle ne permit donc pas aux lecteurs chinois de saisir la démarche hypothético-déductive de la géométrie grecque due à Euclide[7]. En revanche, les néologismes chinois formés pour décrire les figures géométriques s'imposèrent durablement. En 1613, fut publiée une présentation du calcul écrit dont la première partie était une adaptation de l'*Epitome arithmeticae practicae* de Clavius, les chiffres de l'original étant transcrits avec des caractères chinois ce qui donne une

---

[5] Voir Landry-Deron (2013) sur Matteo Ricci dans la Chine des Ming.

[6] Voir Spence (1986) sur cet épisode.

[7] Voir Matzoff (1987: 99-108) et Jami (1998). En 1773, plus d'un siècle et demi après cette publication, l'empereur Qianlong (1711-1799) décida malgré tout que cette traduction des *Éléments* intègrerait la bibliothèque impériale « *Siku quanshu* ». C'est qu'entre temps de nouveaux venus avaient apporté une meilleure compréhension du texte. En effet, en 1685 le roi Louis XIV (1638 –1715) manda en Chine ceux que l'on devait appeler « les mathématiciens du roi » (Landry-Deron 2001) : par ordre de séniorité Jean de Fontaney (1643–1710) qui était régent de mathématiques au Collège Louis-le-Grand, Guy Tachard (1651–1712), Jean-François Gerbillon (1654–1707), Louis Le Comte (1655–1728), Joachim Bouvet (1656–1730) et Claude de Visdelou (1656–1737). Tachard resta au Siam et les cinq autres jésuites français entrèrent à Pékin en 1688. Sur le plan politique, le patronage portugais de la mission jésuite en Chine s'en trouvait ébranlé ; dans le domaine scientifique, les mathématiciens français arrivaient en héritier du père Ignace Gaston Pardies (1636-1673) qui avait précédé Jean de Fontaney à la chaire de mathématiques de Louis-le-Grand et avait écrit une adaptation des *Éléments* plus pédagogique que la traduction de Clavius.

représentation numérique décimale et positionnelle similaire à celle en chiffres indo-arabes.

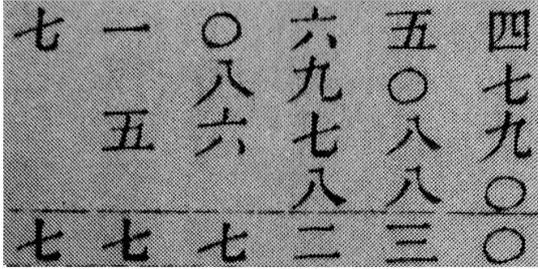

Au début du XVIII[e] siècle, les jésuites accordaient aux chrétiens chinois le droit de célébrer les rites civils dédiés à Confucius et à l'empereur, mais le pape condamna cette pratique. Cette attitude de la papauté conduisit le pouvoir impérial à limiter les droits initialement accordés aux missionnaires dont l'influence déclina[8]. On l'a vu, ils avaient apporté une contribution majeure dans le domaine calendaire. L'introduction des savoirs occidentaux suscita un regain d'intérêt pour les mathématiques de la part des Chinois et les poussa même à relire les textes chinois anciens tombés en désuétude.

Dans les années 1850, Li Shanlan (1811–1882) introduisit les calculs algébrique et différentiel en adaptant des manuels britanniques avec les missionnaires protestants Alexander Wylie (1815-1887) et Joseph Edkins (1823-1905). Il inventa ses propres symboles pour noter les divers opérateurs. Les inconnues et les indéterminées étaient désignées par des caractères en reprenant l'usage des algébristes chinois du XIII[e] siècle. Les nombres étaient écrits à l'aide de notations chinoises. Dans les années 1870, Hua Hengfang (1833–1902), qui traduisait avec John Fryer (1839-1928), continuait d'employer les notations de Li Shanlan. En revanche, à la même époque au Japon, les réformateurs du début de l'ère Meiji (1868-1912) adoptèrent les notations occidentales, et aussi l'essentiel des néologismes forgés par les traducteurs chinois[9]. Les Japonais estimaient que s'approprier le système de notation déjà disponible serait plus rapide que d'en construire un nouveau. En Chine, le même argument finit par l'emporter dans l'effervescence du « mouvement du 4 mai 1919 » qui visait à moderniser la Chine.

---

[8] Sur cette « affaire des rites chinois », on se réfère à Étiemble (1988: vol. 1, 80-293).

[9] Les missionnaires étrangers étaient bannis du Japon où les lettrés se chargeaient eux-mêmes de traduire les ouvrages européens. À l'ère Meiji (1868-1912), le rythme des traductions s'accéléra. S'ils disposaient d'une version chinoise, les traducteurs japonais s'en servaient et s'approprièrent ainsi une partie du lexique scientifique chinois constitué dans les traductions des jésuites (Bréard 2001, 2004) : des néologismes écrits en caractères chinois passaient directement du chinois au japonais et un nouveau vocabulaire scientifique en partie commun aux deux langues se constituait. Au XX[e] siècle, les étudiants chinois au Japon ne percevaient pas toujours si un terme scientifique japonais était un authentique néologisme japonais ou s'il avait été initialement conçu en Chine.

Dès lors les mathématiciens chinois étaient entrés dans le dialogue scientifique mondialisé des temps modernes. Et quand Li Yan (1892-1963) et Qian Baozong (1892-1974) écrivirent dans les années 1920-1950 sur l'histoire des mathématiques en Chine, ils présentèrent les résultats anciens avec les nouvelles notations. Les jeunes générations, qui avait suivi ses études en Chine ou qui revenaient du Japon, des États-Unis ou d'Europe, n'en connaissaient pas d'autres.

Sur l'illustration suivante on trouve des documents chinois et français d'aujourd'hui. Les figures sont identiques, on y voit un repère où est tracée une parabole dont le sommet est désigné par une lettre de l'alphabet latin. Les nombres sont écrits dans la forme européenne des chiffres indo-arabes. Le même symbolisme est utilisé pour l'expression algébrique de la fonction du second degré.

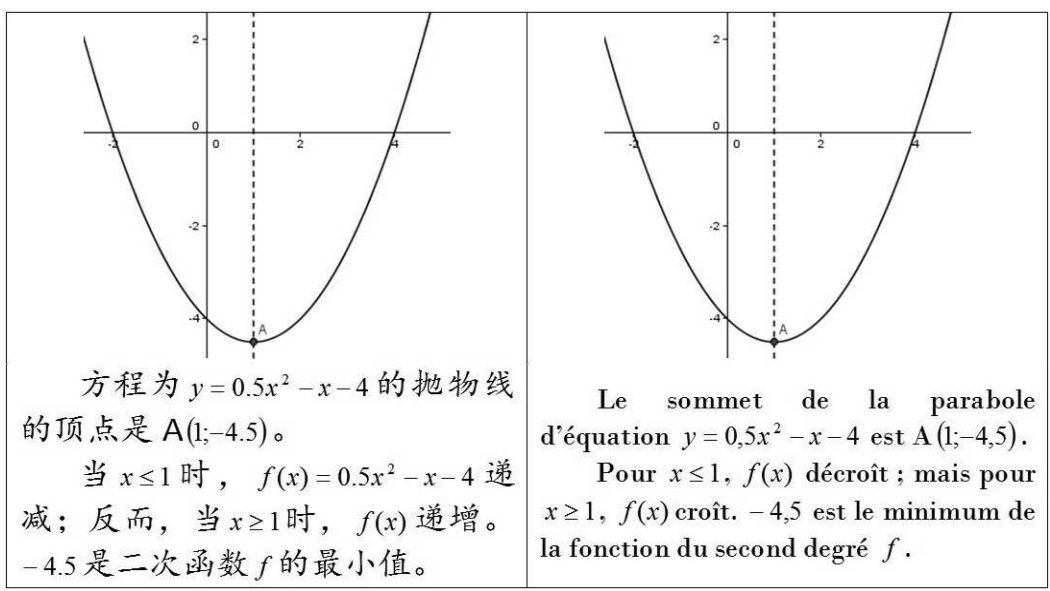

方程为 $y = 0.5x^2 - x - 4$ 的抛物线的顶点是 $A(1;-4.5)$。

当 $x \leq 1$ 时，$f(x) = 0.5x^2 - x - 4$ 递减；反而，当 $x \geq 1$ 时，$f(x)$ 递增。$-4.5$ 是二次函数 $f$ 的最小值。

Le sommet de la parabole d'équation $y = 0,5x^2 - x - 4$ est $A(1;-4,5)$.

Pour $x \leq 1$, $f(x)$ décroît ; mais pour $x \geq 1$, $f(x)$ croît. $-4,5$ est le minimum de la fonction du second degré $f$.